\documentclass[reqno]{amsart}

\usepackage{geometry}
\usepackage[english]{babel}
\usepackage{latexsym, mathrsfs, verbatim}
\usepackage{amsfonts, amsthm, amsmath, amssymb, bm}
\usepackage{hyperref}
\usepackage{amsaddr}
\usepackage{rotating}

\newtheorem{mainthm}{Theorem}

\newtheorem{prop}{Proposition}[section]
\newtheorem{lem}[prop]{Lemma}

\newtheorem{thm}[prop]{Theorem}

\theoremstyle{definition}

\newtheorem{defi}[prop]{Definition}
\newtheorem{ex}[prop]{Example}
\newtheorem{prob}{Problem}

\def\Z{\mathbb{Z}}

\def\equad{\quad \text{ and } \quad}

\def\H{\mathrm{H}}
\def\IHS{\mathrm{IHS}}
\def\supp{\mathrm{supp}}
\def\skel{\mathrm{skel}}
\def\T{\mathsf{t}}

\def\dcup{\sqcup}

\numberwithin{equation}{section}

\usepackage[table,dvipsnames]{xcolor}
 \usepackage{datetime}
 
 \def\sha{\cellcolor{black!10}}

\hypersetup{ 
 colorlinks =true,
 linkcolor=Sepia,
 citecolor=Plum,
 urlcolor=olive,
 }

\begin{document}

\title{An update on the existence of integer Heffter arrays}

\author[Fiorenza Morini]{Fiorenza Morini}
\address{Dipartimento di Scienze Matematiche, Fisiche e Informatiche, Universit\`a di Parma,\\
Parco Area delle Scienze 53/A, 43124 Parma, Italy}
\email{fiorenza.morini@unipr.it}

\author[Marco Antonio Pellegrini]{Marco Antonio Pellegrini}
\address{Dipartimento di Matematica e Fisica, Universit\`a Cattolica del Sacro Cuore,\\
Via della Garzetta 48, 25133 Brescia, Italy}
\email{marcoantonio.pellegrini@unicatt.it}

\begin{abstract}
An integer Heffter array $\H(m,n;s;k)$ is
an $m\times n$ partially filled array whose entries are the elements of a subset $\Omega\subset \Z$ such that $\{\Omega,-\Omega\}$ is a partition of the set $\{1,2,\ldots,2nk\}$ and such that the following conditions are satisfied: each row contains $s$ filled cells, each column contains $k$ filled cells, 
the elements in every row and column  add up to $0$. 
It was conjectured by Dan Archdeacon  that 
an integer $\H(m,n;s;k)$ exists if and only if
$ms=nk$, $3\leqslant s \leqslant n$, $3\leqslant k\leqslant m$ and $nk\equiv 0,3\pmod 4$. 
In this paper, we provide new constructions of these objects that allow us to prove the validity of  Archdeacon's conjecture in each admissible case, except when $k=3,5$ and 
$s\not \equiv 0\pmod 4$ is such that 
$\gcd(s,k)=1$.
\end{abstract}

\keywords{Integer Heffter array; Heffter array set; zero-sum block}
\subjclass{05B20, 05B30}

\maketitle

\section{Introduction}

In 2015 Dan Archdeacon introduced the concept of a Heffter array. 
We recall, in particular, that an integer Heffter array $\H(m,n;s;k)$ is
an $m\times n$ partially filled array whose entries are the elements of a subset $\Omega\subset \Z$ such that $\{\Omega,-\Omega\}$ is a partition of the set $\{1,2,\ldots,2nk\}$ and such that the following conditions are satisfied: each row contains $s$ filled cells, each column contains $k$ filled cells, 
the elements in every row and column  add up to $0$. 
It is easy to see that necessary conditions for the existence of these objects are:
$ms=nk$, $3\leqslant s \leqslant n$, $3\leqslant k\leqslant m$ and $nk\equiv 0,3\pmod 4$. Archdeacon conjectured in \cite{A} that these conditions are also sufficient.
    
One of the main motivations to study the existence of an (integer) $\H(m,n;s,k)$ is that, under some special conditions, such arrays allow one  to construct an orientable embedding of the complete graph $K_{2nk+1}$  in which every edge lies on a face of size $s$ and a face of size $k$ \cite[Theorem 1.1]{A}.
    
The first papers dealing with the existence of 
 integer Heffter arrays considered the tight case (i.e., when $m=k$ and $n=s$)  and the square case (i.e., when $m=n$ and $s=k$),
see \cite{ABD} and \cite{ADDY,DW}, respectively. 
In \cite{MP} the authors solved the case when the parameters $s$ and $k$ are both even.
So, for the remaining cases, one can clearly assume that $k$ is odd.
The same authors also proved the following two results.

\begin{thm}\cite[Theorem 1.4]{MP3}\label{diago}
Let $m, n, s, k$ be four integers such that $3 \leqslant s \leqslant n$, 
$3 \leqslant k \leqslant m$ and $ms = nk$. 
Set $d=\gcd(s,k)$. There exists an integer Heffter array $\H(m, n; s, k)$
in each of the following cases:
\begin{itemize}
    \item[$(1)$] $d\equiv 0 \pmod 4$;
    \item[$(2)$] $d\equiv 1 \pmod 4$ with $d\geqslant 5$ and $nk\equiv 3 \pmod 4$;
    \item[$(3)$] $d\equiv 2 \pmod 4$ and $nk\equiv 0 \pmod 4$;
    \item[$(4)$] $d\equiv 3 \pmod 4$ and 
    $nk\equiv 0,3\pmod 4$.
\end{itemize}
\end{thm}

\begin{thm}\cite[Theorem 1.5]{MP3}\label{so}
Let $m, n, s, k$ be four integers such that $3 \leqslant s \leqslant n$, 
$3 \leqslant k \leqslant m$ and $ms = nk$. 
If $s \equiv  0 \pmod{4}$ and $k \neq 5$ is odd, then
there exists an integer Heffter array $\H(m, n; s, k)$. 
\end{thm}

Theorem \ref{diago} left open the case when 
$nk\equiv 0\pmod 4$ and $d\geq 5$ is such $d\equiv 1 \pmod 4$.
The case when $\gcd(s,k)=1$ was considered  in \cite{PT}.

\begin{thm}\cite[Corollary 1.6]{PT}\label{cop}
Let $m, n, s, k$ be four integers such that 
$3 \leqslant s \leqslant n$, $3 \leqslant k \leqslant m$,
$ms = nk$ and $nk \equiv 0, 3 \pmod 4$. 
There exists an integer $\H( m, n; s, k )$ whenever 
$s, k$ are coprime integers such that $k \geqslant 7$ is odd and
$s \neq 3, 5, 6, 10$. 
\end{thm}

The previous two results were based on the existence of particular sets of blocks.
Hence, inspired also by the ideas of \cite{C,F13} where ($\Gamma$-)magic rectangle sets were considered (see also \cite{MP4}), we give the following.

\begin{defi}\label{IHS}
An \emph{integer Heffter array set} $\IHS(m,n; c)$ is a collection of $c$ arrays of size $m\times n$  such that
\begin{itemize}
\item[{\rm (a)}] the entries are the elements of a subset $\Omega\subset \Z$  such that $\{\Omega, -\Omega\}$ is a partition of $\{1,\ldots,2mnc\}$;
\item[{\rm (b)}] every $\omega \in \Omega$ appears  once and in a unique array;
\item[{\rm (c)}] for every array, the sum of the elements in each row  and in each column is $0$.
\end{itemize}
\end{defi}

Also in this case, it is clear that the following conditions are necessary for the existence of
an $\IHS(m,n;c)$:
$m,n\geqslant 3$ and $mnc\equiv 0,3 \pmod 4$.

The main goal of this paper is to close some of the cases left open by the previous theorems. So, we first show that Archdeacon's conjecture holds when 
$\gcd(s,k)>1$.

\begin{mainthm}\label{mainA}
Let $m, n, s, k$ be four integers such that $3 \leqslant s \leqslant n$, 
$3 \leqslant k \leqslant m$ and $ms = nk$.
Let $d=\gcd(s,k)$.
If $d\equiv 1 \pmod 4$, $d\geqslant 5$ and $nk\equiv 0 \pmod 4$, then there exists an integer Heffter array $\H(m, n; s, k)$.
\end{mainthm}

Next, we  construct integer Heffter array sets to prove the following.

\begin{mainthm}\label{mainB}
Let $m,n,c$ be three positive integers such that $m,n\geqslant 3$ and $mnc\equiv 0,3 \pmod 4$.
There exists an $\IHS(m,n;c)$ in each of the following cases:
\begin{itemize}
\item[$(1)$] $m,n$ are even integers;
\item[$(2)$] $m\equiv 0 \pmod 4$ and $n$ is an odd integer;
\item[$(3)$] $m\equiv 2 \pmod 4$ and $n\geqslant 7$ is an odd integer; 
\item[$(4)$] $m,n\geqslant 7$ are odd integers.
\end{itemize}
\end{mainthm}

To conclude, we remark that the validity of Archdeacon's conjecture on the existence of an integer Heffter array $\H(m,n;s,k)$ is a problem which remains open only
when $k=3,5$ and $s\not \equiv 0 \pmod 4$ is such that $\gcd(s,k)=1$.
Indeed, as a consequence of the previous results in  \cite{MP,MP3,PT} and of Theorems \ref{mainA} and \ref{mainB} we have the following.

\begin{mainthm}\label{mainC}
Let $m, n, s, k$ be four integers such that $3 \leqslant s \leqslant n$, 
$3 \leqslant k \leqslant m$, $ms = nk$ and $nk\equiv 0,3 \pmod 4$. There exists an integer Heffter array $\H(m, n; s, k)$ in each of the following cases:
\begin{itemize}
\item[$(1)$] $\gcd(s,k) \neq 1$;
\item[$(2)$] $s\equiv 0 \pmod 4$;
\item[$(3)$] $s, k$ are coprime integers such that $k \geqslant 7$ is odd and $s \neq 3, 5$. 
\end{itemize}
\end{mainthm}

A solution to the following two problems would allow to provide a complete proof of Archdeacon's conjecture on the existence of integer Heffter arrays.

\begin{prob}
Construct an $\IHS(m,3;c)$ for any $c\geqslant 1$ and any $m\geqslant 3$ such that $m$ is odd, $\gcd(m,3)=1$ and $mc\equiv 0,1 \pmod 4$.
\end{prob}

\begin{prob}
Construct an $\IHS(m,5;c)$ for any $c\geqslant 1$ and any $m\geqslant 5$ such that $m$ is odd, $\gcd(m,5)=1$ and $mc\equiv 0,3 \pmod 4$.
\end{prob}

\section{Notation}

Given an integer $q\geqslant 1$, if $a,b$ are two integers such that $a\equiv b \pmod{q}$, then we use the notation
$$[a,b]_q =\left\{a+iq\mid 0\leqslant i \leqslant \frac{b-a}{q}\right\},$$ whenever 
$a\leqslant b$. If $a>b$, then $[a,b]_q=\varnothing$.
If $q=1$, we simply write $[a,b]$.
For $q\in\{1,2\}$, any $\ell$-subset of $\mathbb{Z}$ of the form $[x, x+(\ell-1)q]_{q}$ will be called an $\ell$-set of type $q$.

Let $A$ be a partially filled array with  integer
entries. The support of $A$, denoted by $\supp(A)$, is defined as the list of the absolute values
of the entries of $A$.
Given a set $\mathfrak{S} = \{ A_1, A_2 ,
\ldots, A_r\}$ of partially filled arrays with  integer entries, we set 
$\supp ( \mathfrak{S} ) = \cup_i \supp (A_i )$.
The skeleton of $A$, denoted by $\skel(A)$, is
the set of its filled  cells.

We denote by $\sigma_r(A)$ and $\sigma_c(A)$, respectively, the sequences of the sums of the elements of each row and of each  column of  $A$. In particular, if $A$ has no empty cells and $\sigma_r(A)$ and $\sigma_c(A)$ are sequences of zeroes, we say that $A$ is a zero-sum block.

An integer $\H(m,n;s,k)$ is said to be shiftable if every row and every column  contains the same number of positive and negative entries.
We recall that a shiftable $\H(m,n;s,k)$ 
exists if and only if the following conditions are satisfied: $4\leqslant s \leqslant n$, $4\leqslant k\leqslant m$, $ms=nk$
and $s\equiv k \equiv 0\pmod 2$, see
\cite{MP}.
Note that, given a shiftable $\H(m,n;s,k)$, say $A$,
and a positive integer $\alpha$, we can replace every positive entry $x$ with $x+\alpha$ and every negative entry $-y$ with $-y-\alpha$, obtaining a partially filled array, denoted by $A\pm \alpha$, satisfying the definition of integer Heffter array, except for its support. In fact, $\skel(A)=\skel(A\pm \alpha)$ and $\supp(A\pm \alpha)=[\alpha+1, \alpha+nk]$.

If $A$ and $B$ are two partially filled arrays of the same size  such that $\skel(A)\cap \skel(B)=\varnothing$, we denote by $A\oplus B$ the partially filled array obtained by overlapping $A$ and $B$.

To simplify  our notation, we will write $\H(n;k)$ instead of $\H(n,n;k,k)$,
and $\H(m,n)$ instead of $\H(m,n;n,m)$.

\section{Proof of Theorem \ref{mainA}}

In this section, we prove the existence of an integer Heffter array $\H(m,n;s,k)$ when $d=\gcd(s,k)$ is such that $d\geqslant 5$ and $d\equiv 1 \pmod 4$. Thus, we can write 
\begin{equation}\label{s1k1}
m=ek_1, \quad n=es_1, \quad s=ds_1 \equad 
k=dk_1
\end{equation}
for some $e\geqslant d$, where $s_1,k_1\geqslant 1$ and $\gcd(s_1,k_1)=1$.
Since the square case has already been solved, we can assume $e>d$.
We start by considering the case $d=5$.

\begin{lem}\label{d5}
Suppose that $k_1,s_1$ are two positive and coprime integers. For every $e>5$ such that $es_1k_1\equiv 0 \pmod 4$, there exists an 
integer $\H(ek_1, es_1; 5s_1, 5k_1)$.
\end{lem}

\begin{proof}
Set $4N=es_1k_1$. Let $A$ be an integer $\H(4N; 5)$ as constructed in \cite[Theorem 4.2]{DW}.
Then
$$\skel(A)=\left\{(i,j): 1\leqslant i,j\leqslant 4N,\; j-i \in \mathcal{A}\right\},$$
where $\mathcal{A}=\{-(4N-1), -(4N-2), -2N, -(2N-1), 0,1,2,  2N, 2N+1\}$.
Note that the difference between any two integers of $\mathcal{A}$ is an element of 
$$\begin{array}{rcl}
\mathcal{B} & = &  \pm    \{ 0, 1, 2, 2N-2, 2N-1,  2N, 2N+1, 2N+2, 4N-2, 4N-1,  \\
&& 4N,  4N+1,  6N-2, 6N-1,  6N \}.
  \end{array}$$
As $e>5$, the only multiples of $e$ in $\mathcal{B}$ belong to the subset $\pm \{0, 2N, 4N, 6N\}$.

Since $s_1$ and $k_1$ cannot be both even, we may assume that $k_1\geqslant 1$ is odd.
Define the function
$$\pi: \skel(A) \to \{(u,v):  1\leqslant u\leqslant ek_1, \; 1\leqslant v\leqslant es_1 \}$$
as follows. Given $(i,j)\in \skel(A)$, there exist four uniquely determined integers $q_i,q_j,u_i,v_j$ such that
$$1\leqslant u_i\leqslant ek_1,\quad 1\leqslant v_j\leqslant es_1, \quad i+z_i = q_i(ek_1)+ u_i \equad j = q_j (es_1)+v_j,$$
where
$$z_i=\left\{\begin{array}{ll}
3 & \text{ if } s_1 \text{ is even and } 2N+1\leqslant i \leqslant4N, \\
0 & \text{ otherwise}.
             \end{array}\right.$$
Set $\pi(i,j)=(u_i,v_j)$. We show that $\pi$ is an injective function.
Suppose that $(i,j)$ and $(a,b)$ are two elements of $\skel(A)$ such that $\pi(i,j)=\pi(a,b)$.
Using the previous notation, we have $i+z_i=a+z_a +(q_i-q_a)(e k_1)$ and $j = b+(q_j-q_b)(e s_1)$.
Writing $x=q_i-q_a$ and $y=q_j-q_b$, this implies that  $j-i+(z_a-z_i) = b-a + e (ys_1 - xk_1)$,
where $0\leqslant x < s_1$ and $0\leqslant y < k_1$.
In particular, $(j-i)-(b-a)= e (ys_1 - xk_1)-(z_a-z_i)$ is an element of $\mathcal{B}$.

First, suppose that $s_1$ is odd. In this case, $z_i=z_a=0$.
So $(j-i)-(b-a)= e (ys_1 - xk_1)$ is an element of $\mathcal{B}$ which is a 
multiple of $e$.
Since the product $s_1k_1$ is odd, the integers $2N$ and $6N$ are not multiples of $e$:
by the above, we obtain that either $ys_1 - xk_1=0$ or $ys_1 - xk_1=\pm (s_1k_1)$.
As $\gcd(s_1,k_1)=1$, in  both cases we get that $s_1\mid x $ and $k_1\mid y$, whence $x=y=0$.
We conclude that $i=a$ and $b=j$, proving the injectivity of $\pi$.

Now, suppose that $s_1$ is even. If $1\leqslant i,a\leqslant 2N$ or $2N+1\leqslant i,a\leqslant 4N$, then 
$z_i=z_a$ and we can proceed as before, obtaining that one of the following cases holds:
\begin{itemize}
\item[$(1)$] $ys_1 - xk_1=0$;
\item[$(2)$] $ys_1 - xk_1=\pm \frac{s_1}{2}k_1$;
\item[$(3)$] $ys_1 - xk_1=\pm s_1k_1$;
\item[$(4)$] $ys_1 - xk_1=\pm \frac{3s_1}{2} k_1$.
\end{itemize}
In all four cases, we obtain that  $k_1$ divides $y$  since $\gcd(s_1,k_1)=1$, whence $y=0$.
In case $(1)$ we easily obtain $x=0$; 
cases $(2)$, $(3)$ and $(4)$ are excluded  because  the hypotheses on $i$ and $a$ imply that $0\leqslant x < \frac{s_1}{2}$.
Finally, assume that $1\leqslant i\leqslant 2N$ and $2N+1\leqslant a\leqslant 4N$.
Then, $(j-i)-(b-a)+3$  
is a 
multiple of $e$. It is clear that this cannot happen.
This proves the injectivity of $\pi$.

Finally, we can construct an $(ek_1)\times (es_1)$ array $H$ such that $\skel(H)=\pi(\skel(A))$, filling the cell $(u,v)$ with the entry of the cell $(i,j)$ of $A$, where $\pi(i,j)=(u,v)$.
Thus, the resulting array $H$ is an integer Heffter array $\H(ek_1,es_1;5s_1,5k_1)$. Indeed, $\supp(A)=\supp(H)=[1,5es_1k_1]$. Moreover, each row of $H$ is 
obtained by taking $s_1$ rows of the array $A$, and hence contains $5s_1$ filled cells whose entries add up to zero. Analogously, each column of  $H$ is obtained by taking $k_1$ columns of $A$, and hence  contains $5k_1$ filled cells with zero-sum entries.
\end{proof}

\begin{proof}[Proof of Theorem \ref{mainA}]
Keeping the notation \eqref{s1k1},  by Lemma \ref{d5} it suffices to show  how to construct an integer $\H(ek_1, es_1; ds_1, dk_1)$, when $e> d > 5$ and $es_1k_1 \equiv 0 \pmod 4$. So, write $d = 5+4t$ where $t>0$.

By easily adapting the construction described in the proof of \cite[Proposition 3.5]{MP}, we construct a shiftable 
$\H(ek_1,es_1; 4ts_1, 4tk_1)$, say $C$, whose skeleton is 
$$\{(a,b): 1\leqslant a\leqslant ek_1,\;
1\leqslant b\leqslant es_1, \; 
b-a \equiv \ell \pmod e, \; \ell \in \mathcal{L}\},$$
where 
$$\mathcal{L}=\left\{\begin{array}{ll}
{}[3,2+2t] \dcup [e-2t,e-1] & \text{if $s_1$ is odd},\\
{}[3,4t+2] & \text{if $s_1$ is even}.
\end{array}\right.$$
Furthermore, let $H$ be the integer $\H(ek_1,es_1;5s_1,5k_1)$
constructed in the proof of Lemma \ref{d5}.
Keeping the previous notation, let $(u_i,v_j)=\pi(i,j)\in\skel(H)$.

Suppose that $s_1$ is odd. 
Since $v_j-u_i \equiv j-i \pmod e $ and $j-i \in \mathcal{A}$, it follows that $v_j-u_i\equiv \mu \pmod e$, where  $\mu \in \left\{0,1,2, \frac{e}{2},\frac{e}{2}+1\right\}$. In fact, $2N=\frac{e}{2}s_1k_1\equiv \frac{e}{2} \pmod e$, because $s_1k_1$ is odd, and $4N=e(s_1 k_1)$.
Since $e>4t+5$, we have  $2< 3 < 2+2t <\frac{e}{2}<\frac{e}{2}+1< e-2t$, proving that $\skel(H)\cap \skel(C)=\varnothing$.

Suppose now that $s_1$ is even. 
Since $v_j-u_i \equiv j-i-z_i \pmod e $,
it follows that  $v_j-u_i\equiv \eta \pmod e$, where  
$\eta \in \{0,1,2, e-3,e-2,e-1\}$.
In fact, $2N=e \frac{s_1}{2} k_1$.
As $e>4t+5$, we have  $2< 3 < 4t+2 <e-3$, proving 
also in this case that $\skel(H)\cap \skel(C)=\varnothing$.

Taking $H \oplus (C\pm 5ek_1s_1)$ we obtain an integer Heffter array $\H(m,n;s,k)$.
\end{proof}

\begin{ex} 
Figure \ref{fig:H18} shows an integer $\H(20,10;9,18)$  obtained by following the proof of Theorem \ref{mainA}. 
Specifically, we have constructed  an integer 
$\H(10,20;10,5)$, say $H$, according to Lemma \ref{d5}, together with a shiftable array $\H(10,20;8,4)$, say $C$. 
The array shown in Figure \ref{fig:H18} is the transpose of  $H\oplus (C\pm 100)$, with the filled cells of $C\pm 100$ highlighted in grey. \end{ex}

\begin{figure}[ht]
$\begin{array}{|c|c|c|c|c|c|c|c|c|c|}\hline
 58 & -40 & 1 & 19 & \sha -177 &\sha  \sha 178 & \sha 101 & \sha -102 &  & -38 \\\hline
 -18 & -79 & 96 & -20 & 21 &\sha  -179 &\sha  180 & \sha 103 &\sha  -104 &  \\\hline
 78 & -17 & -62 &  & -22 & 23 &\sha  -141 &\sha  142 & \sha 105 &\sha  -106 \\\hline 
 \sha -108 & 95 & -16 & -80 &  & -24 & 25 & \sha -143 &\sha  \sha 144 &\sha  107 \\\hline
 \sha 109 &\sha  -110 & 77 & -15 & -63 &  & -26 & 27 & \sha -145 &\sha  146 \\\hline
\sha 148 & \sha 111 &\sha  -112 & 94 & -14 & -81 &  & -28 & 29 &\sha  -147 \\\hline
\sha  -149 &\sha  150 &\sha  113 & \sha -114 & 76 & -13 & -64 &  & -30 & 31 \\\hline
 33 &\sha  -151 & \sha 152 &\sha  115 &\sha  -116 & 93 & -12 & -82 & & -32\\\hline
-34 & 35 & \sha -153 & \sha 154 &\sha  117 &\sha  -118 & 75 & -11 & -65 & \\\hline
  & 100 & -99 &\sha  -155 &\sha  156 & \sha 119 &\sha  -120 & 92 & -10 & -83 \\\hline
-61 &  & -37 & -73 &\sha  -157 & \sha 158 &\sha  121 &\sha  -122 & 74 & 97 \\\hline
-57 & 56 &  & 8 & -91 & \sha -159 &\sha  160 &\sha  123 &\sha  -124 & 84 \\\hline 
  & -55 & 54 & 66 & 7 & -72 &\sha  -161 &\sha  162 &\sha  125 &\sha  -126 \\\hline 
\sha -128 &  & -53 & 52 & 85 & 6 & -90 &\sha  -163 &\sha  164 &\sha  127 \\\hline 
 \sha 129 &\sha  -130 &  & -51 & 50 & 67 & 5 & -71 &\sha  -165 &\sha  166 \\\hline
 \sha 168 &\sha  131 &\sha  -132 & & -49 & 48 & 86 & 4 & -89 &\sha  -167 \\\hline
 \sha -169 & \sha 170 &\sha  133 &\sha  -134 &  & -47 & 46 & 68 & 3 & -70 \\\hline
 -88 & \sha -171 &\sha  172 & \sha 135 &\sha  -136 &  & -45 & 44 & 87 & 2 \\\hline
 -9 & -59 &\sha  -173 & \sha 174 &\sha  137 &\sha  -138 & & -43 & 42 & 69 \\\hline
 98 & -36 & 39 &\sha  -175 & \sha 176 &\sha  139 &\sha  -140 &  & -41 & -60 \\\hline
\end{array}$
\caption{An integer $\H(20,10;9,18)$.}\label{fig:H18}
\end{figure}

\section{Proof of Theorem \ref{mainB}}

In this section we consider the existence of some integer Heffter array sets.

\begin{lem}\label{45even}
There exists an $\IHS(4,5;c)$ for every even $c\geqslant 2$.
\end{lem}

\begin{proof}
 Let $c=2t+2$ where $t\geqslant 0$.
For every $i \in [0,t]$, define
$$\begin{array}{rcl}
A_{2i} & =& \begin{array}{|c|c|c|c|}\hline
1+16i &  -(3+16i) & -(5+16i) & 7+16i \\ \hline
-(2+16i) & 4+16i & 10+16i & -(12+16i) \\\hline
6+16i & -(8+16i) & -(14+16i) & 16+16i \\ \hline
24t+24-8i & -(24t+23-8i)  & -(24t+22-8i) & 24t+21-8i \\ \hline 
-(24t+29+8i) & 24t+30+8i & 24t+31+8i & -(24t+32+8i) \\ \hline
\end{array},\\ \\[-9pt]
A_{2i+1} & =& \begin{array}{|c|c|c|c|}\hline
24t+20-8i & -(24t+19-8i) & -(24t+18-8i) & 24t+17-8i \\ \hline
-(24t+25+8i) & 24t+26+8i & 24t+27+8i & -(24t+28+8i) \\ \hline
32t+33+8i & -(32t+34+8i) & -(32t+35+8i) & 32t+36+8i \\ \hline
-(32t+37+8i) & 32t+38+8i & 32t+39+8i & -(32t+40+8i) \\ \hline
9+16i & -(11+16i) & -(13+16i) & 15+16i \\ \hline
\end{array}.
\end{array}$$

The set $\mathcal{A}=\{A_{2i}^\T,A_{2i+1}^\T: i \in [0,t]\}$ consists of $2t+2$ zero-sum blocks of size $4\times 5$ and has support equal to $[1, 40t+40]$. In fact, for every $i$, we have  $\supp(A_{2i})\dcup \supp(A_{2i+1})=
[1+16i, 16+16i] \dcup [24t+17-8i ,24t+24-8i ]
\dcup [24t+25+8i, 24t+32+8i]\dcup [32t+33+8i, 32t+40+8i]$.
We conclude that $\mathcal{A}$ is an $\IHS(4,5;2t+2)$.
\end{proof}

\begin{lem}\label{45odd}
There exists an $\IHS(4,5;c)$ for every odd $c\geqslant 1$.
\end{lem}

\begin{proof}
Let $c=2t+1$ where $t\geqslant 0$.
Define
$$A = \begin{array}{|c|c|c|c|}\hline
16t+1 & -(16t+3) & -(16t+5) & 16t+7 \\ \hline
-(16t+10)& 16t+9 & 16t+12 & -(16t+11) \\ \hline
-(40t+17) & 40t+16 & -(40t+18) & 40t+19 \\ \hline
40t+20 & -(40t+14) & -4 & -2 \\ \hline 
6 & -8 & 40t+15 & -(40t+13) \\ \hline
\end{array}.$$
Furthermore, for every $i \in [0, t-1]$, define
$$\begin{array}{rcl}
B_{2i} & =& \begin{array}{|c|c|c|c|}\hline
1+16i & -(3+16i) & -(5+16i)& 7+16i \\ \hline
-(10+16i) & 12+16i & 18+16i & -(20+16i) \\ \hline
14+16i & -(16+16i) & -(22+16i) & 24+16i  \\ \hline
24t+12-8i & -(24t+11-8i) & -(24t+10-8i) & 24t+9-8i \\ \hline
-(24t+17+8i)& 24t+18+8i & 24t+19+8i & -(24t+20+8i) \\\hline
\end{array},\\ \\[-9pt]
B_{2i+1} & =& \begin{array}{|c|c|c|c|}\hline
24t+8-8i & -(24t+7-8i) & -(24t+6-8i) & 24t+5-8i \\ \hline
-(24t+13+8i) & 24t+14+8i & 24t+15+8i & -(24t+16+8i ) \\ \hline
32t+13+8i &  -(32t+14+8i) & -(32t+15+8i) & 32t+16+8i \\ \hline
-(32t+17+8i) & 32t+18+8i & 32t+19+8i & -(32t+20+8i)  \\ \hline
9+16i & -(11+16i) & -(13+16i) & 15+16i \\ \hline
\end{array}.
\end{array}$$

The set $\mathcal{A}=\{A^\T\}\cup \{B_{2i}^\T, B_{2i+1}^\T: i \in [0,t-1]\}$ consists of $2t+1$ zero-sum blocks of size $4\times 5$ and has support equal to $
[1,  40t+20]$.
In fact, $\supp(A)=[2,8]_2 \dcup [16t+1,16t+7 ]_2
\dcup [16t+9, 16t+12]\dcup [40t+13, 40t+20]$.
Furthermore, for every $i$, we have  $\supp(B_{2i})\dcup \supp(B_{2i+1})=
[1+16i, 15+16i ]_2 \dcup [10+16i, 24+16i  ]_2 \dcup [24t+ 5-8i, 24t+12-8i]
\dcup [24t+13+8i, 24t+20+8i]\dcup [32t+13+8i, 32t+20+8i ]$.
We conclude that $\mathcal{A}$ is an $\IHS(4,5;2t+1)$.
\end{proof}

We now consider the existence of an $\IHS(m,n;c)$ when 
$m\equiv 2 \pmod 4$ and $n$ is odd. We start with two auxiliary lemmas.

\begin{lem}\label{BlockA}
Given three positive integers $\alpha,\beta,u$ such that
$\beta\geqslant \alpha+6u-2$,
there exists a set $\mathfrak{A}=\mathfrak{A}(\alpha,\beta,u)$, consisting of  $2\times 3$ matrices such that
\begin{itemize}
\item $|\mathfrak{A}|=u$;
\item $\sigma_r(A)=(0,0)$ and $\sigma_c(A)=(-2,1,1)$ for all $A\in \mathfrak{A}$;
\end{itemize}
and
$$\supp(\mathfrak{A})  = [\alpha,\alpha+4u-2]_2\dcup [\beta-2u+1,\beta] \dcup  [\alpha+\beta, \alpha+\beta+2u-1].$$
\end{lem}

\begin{proof}
For every $j\in [0, u-1]$, define
$$A_j=\begin{array}{|c|c|c|}\hline
\alpha +4j & -( \alpha+\beta+2j) & \beta-2j \\ \hline
-(\alpha+2+4j) & \alpha+\beta+1+2j & -(\beta-1-2j ) \\\hline
\end{array}.$$
The set $\mathfrak{A}=\left\{A_j: j \in [0, u -1]\right\}$ has the required properties.
\end{proof}

\begin{lem}\label{BlockB}
Given two positive integers $\beta$ and $u$ such that
$\beta\geqslant 12u-1 $,
there exists a set $\mathfrak{B}=\mathfrak{B}(\beta,u)$, consisting of  $2\times 3$ matrices, such that
\begin{itemize}
\item $|\mathfrak{B}|=2u$;
\item $\sigma_r(B)=(0,0)$ and $\sigma_c(B)=(-4,2,2)$ for all $B\in \mathfrak{B}$;
\end{itemize}
and
$$\supp(\mathfrak{B})  = [1, 8u-1]_2\dcup [\beta-4u+1,  \beta+4u].$$
\end{lem}

\begin{proof}
For every $j\in [0, u-1]$, define
$$\begin{array}{rcl}
B_{2j} & =& \begin{array}{|c|c|c|}\hline
1 +8j & -( \beta+1+4j) & \beta-4j \\ \hline
-(5+8j) & \beta+3+4j & -(\beta-2-4j ) \\\hline
\end{array},\\[10pt]
B_{2j+1} & =& \begin{array}{|c|c|c|}\hline
3+8j & -( \beta+2+4j) & \beta-1-4j \\ \hline
-(7+8j) & \beta+4+4j & -(\beta-3-4j ) \\\hline
\end{array}.
  \end{array}$$
The set $\mathfrak{B}=\left\{B_{2j}, B_{2j+1}: j \in [0, u -1]\right\}$ has the required properties.
\end{proof}

It will be useful to introduce the following zero-sum blocks $P_1, P_2$ of size $4\times 4$, $Q_1,Q_2,Q_3$ of size $6\times 4$ and $R_1,R_2$ of size $6\times 6$.
So, given the sets 
$\mathcal{X}_i=[x_i+2, x_i+8]_2$, 
$\mathcal{Y}_j=[y_j+1, y_j+4]$, 
$\mathcal{Z}_k=[z_k+1,z_k+8]$
and $\mathcal{W}_{\ell}=[w_\ell+2, w_{\ell}+16]_2$, define:
$$P_1(\mathcal{X}_1,\mathcal{X}_2,\mathcal{X}_3,\mathcal{X}_4)=
\begin{array}{|c|c|c|c|}\hline
  x_1+2 & -(x_1+6) & -(x_2+2) &  x_2+6 \\ \hline
-(x_1+4) &  x_1+8  &  x_2+4   & -(x_2+8) \\\hline
-(x_3+2) & x_3+6  &  x_4+2  & -(x_4+6) \\\hline
x_3+4    & -(x_3+8) & -(x_4+4) & x_4+8\\\hline
\end{array},$$
$$P_2(\mathcal{Y}_1,\mathcal{Y}_2, \mathcal{Y}_3,\mathcal{Y}_4)=\begin{array}{|c|c|c|c|}\hline
  y_1+1 & -(y_1+3) & -(y_2+1) &  y_2+3 \\ \hline
-(y_1+2) &  y_1+4  &  y_2+2   & -(y_2+4) \\\hline
-(y_3+1) & y_3+3  &  y_4+1  & -(y_4+3) \\\hline
y_3+2    & -(y_3+4) & -(y_4+2) & y_4+4\\\hline
\end{array},$$
$$Q_1(\mathcal{X}_1, \mathcal{X}_2, \mathcal{Y}_1,\mathcal{Y}_2,\mathcal{Y}_3,\mathcal{Y}_4)=\begin{array}{|c|c|c|c|}\hline
  x_1+2 & -(x_1+6) & -(x_2+2) &  x_2+6 \\ \hline
-(x_1+4) &  x_1+8  &  x_2+4   & -(x_2+8) \\\hline
-(y_1+1) & y_1+3   & y_2+1 & -(y_2+3) \\\hline
 y_1+2   &-(y_1+4) & -(y_2+2) & y_2+4 \\\hline
 -(y_3+1) & y_3+3 & y_4+1 &  -(y_4+3) \\\hline
 y_3+2 & -(y_3+4) & -(y_4+2) &  y_4+4\\ \hline
      \end{array},$$
$$Q_2(\mathcal{X}_1,\mathcal{Y}_1,\mathcal{Y}_2,\mathcal{Y}_3, \mathcal{Z}_1)=\begin{array}{|c|c|c|c|}\hline
  x_1+2 & -(x_1+4) & -(y_1+1) & y_1+3 \\ \hline
-(x_1+6) &  x_1+8  & y_1+2 & -(y_1+4)\\\hline
-(y_2+1) & y_2+2 & y_3+3 & -(y_3+4) \\\hline
y_2+3 & -(y_2+4) & -(y_3+1) & y_3+2\\\hline
z_1+8    &    z_1+1 & -(z_1+5)   & -(z_1+4)\\ \hline
-(z_1+6) & -(z_1+3) & z_1+2 & z_1+7 \\\hline
      \end{array},$$
$$Q_3( \mathcal{Y}_1,\mathcal{Y}_2,\mathcal{Y}_3,\mathcal{Y}_4, \mathcal{Y}_5, \mathcal{Y}_6)=\begin{array}{|c|c|c|c|}\hline
y_1 +1 & -(y_1+2) & -(y_2+1) & y_2+2 \\ \hline
-(y_1+3) & y_1+4 & y_2+3  & -(y_2+4) \\\hline
-(y_3+1) & y_3+3 & y_4+1 & - (y_4+3) \\\hline
y_3+2 & -(y_3+4) &-(y_4+2) & y_4+4 \\\hline
-(y_5+1) & y_5+3 & y_6+1 & -(y_6+3) \\\hline
y_5+2 & -(y_5+4) & -(y_6+2) & y_6+4\\\hline
\end{array},$$
$$\begin{array}{l}
R_1(\mathcal{X}_1,\mathcal{X}_2,\mathcal{X}_3,\mathcal{W}_1,
\mathcal{W}_2,\mathcal{W}_3)=\\[5pt]
\begin{array}{|c|c|c|c|c|c|}\hline
w_1+2   &  -(w_1+6) & -(w_1+10)  & w_1+12   &  x_1+8 & -(x_1+6) \\\hline
-(w_1+4) &  w_1+14 & w_1+8       & -(w_1+16) & -(x_1+4) & x_1+2 \\\hline
x_2+8 & -(x_2+6) & w_2+2   &  -(w_2+6) & -(w_2+10)  & w_2+12  \\\hline
-(x_2+4) & x_2+2 & -(w_2+4) &  w_2+14 & w_2+8       & -(w_2+16)\\\hline
-(w_3+10)  & w_3+12  & x_3+8 & -(x_3+6)  & w_3+2   &  -(w_3+6)\\\hline
 w_3+8       & -(w_3+16) & -(x_3+4) & x_3+2&  -(w_3+4) &  w_3+14 \\\hline
  \end{array},  \end{array}$$
$$\begin{array}{l}
R_2(\mathcal{Y}_1,\mathcal{Y}_2, \mathcal{Y}_3, \mathcal{Z}_1,
\mathcal{Z}_2,\mathcal{Z}_3)=\\[5pt]
\begin{array}{|c|c|c|c|c|c|}\hline
z_1+1    & -(z_1+3) & -(z_1+5) & z_1+6     &  y_1+4    & -(y_1+3) \\\hline
-(z_1+2) &  z_1+7   & z_1+4    & -(z_1+8)  & -(y_1+2)  & y_1+1 \\\hline
y_2+4    & -(y_2+3) & z_2+1    &  -(z_2+3) & -(z_2+5)  & z_2+6  \\\hline
-(y_2+2) & y_2+1    & -(z_2+2) &  z_2+7    & z_2+4     & -(z_2+8)\\\hline
-(z_3+5) & z_3+6    & y_3+4    & -(y_3+3)  & z_3+1     &  -(z_3+3)\\\hline
 z_3+4   & -(z_3+8) & -(y_3+2) & y_3+1     &  -(z_3+2) &  z_3+7\\ \hline
  \end{array}.  \end{array}$$

\begin{prop}\label{n3}
Let $m,n,c$ be three positive integers such that $m\geqslant 6$ and $n\geqslant 7$.
Suppose that $m\equiv 2 \pmod 4$, $n\equiv 3\pmod 4$  and $c$ is even.
Then, there exists an $\IHS(m,n;c)$.
\end{prop}

\begin{proof}
Let $m=6+4v$, $n=7+4w$ and $c=2t+2$, where $v,w,t\geqslant 0$. 
Our $\IHS(m,n;c)$ will consist of $2t+2$ arrays of type
$$\begin{array}{|c|cc|cc|cc|cccc}\hline
  Y      & \multicolumn{2}{c|}{} & \multicolumn{2}{c|}{} & \multicolumn{2}{c|}{} &\multicolumn{2}{c|}{}  \\ \cline{1-1}
  -X     & \multicolumn{2}{c|}{\smash{\raisebox{0\normalbaselineskip}{$Z_6$}}}  & \multicolumn{2}{c|}{\smash{\raisebox{0\normalbaselineskip}{$Z_6$}}} & \multicolumn{2}{c|}{\smash{\raisebox{0\normalbaselineskip}{$\cdots$}}} & \multicolumn{2}{c|}{\smash{\raisebox{0\normalbaselineskip}{$Z_6$}}}\\  \cline{1-1}
  -X      & \multicolumn{2}{c|}{} & \multicolumn{2}{c|}{} & \multicolumn{2}{c|}{} &\multicolumn{2}{c|}{}   \\ \hline
   X     &  \multicolumn{2}{c|}{} & \multicolumn{2}{c|}{} & \multicolumn{2}{c|}{} &\multicolumn{2}{c|}{}\\ \cline{1-1}
  -X     &   \multicolumn{2}{c|}{\smash{\raisebox{.5\normalbaselineskip}{$Z_4$}}}   &   \multicolumn{2}{c|}{\smash{\raisebox{.5\normalbaselineskip}{$Z_4$}}}   &   \multicolumn{2}{c|}{\smash{\raisebox{.5\normalbaselineskip}{$\cdots$}}}   &   \multicolumn{2}{c|}{\smash{\raisebox{.5\normalbaselineskip}{$Z_4$}}}  \\  \hline 
  \vdots     &     \multicolumn{2}{c|}{\vdots} & \multicolumn{2}{c|}{\vdots} & \multicolumn{2}{c|}{\ddots} &\multicolumn{2}{c|}{\vdots}  \\[1mm] \hline
   X     &  \multicolumn{2}{c|}{} & \multicolumn{2}{c|}{} & \multicolumn{2}{c|}{} &\multicolumn{2}{c|}{}  \\ \cline{1-1}
  -X     &   \multicolumn{2}{c|}{\smash{\raisebox{.5\normalbaselineskip}{$Z_4$}}}   &   \multicolumn{2}{c|}{\smash{\raisebox{.5\normalbaselineskip}{$Z_4$}}}   &   \multicolumn{2}{c|}{\smash{\raisebox{.5\normalbaselineskip}{$\cdots$}}}   &   \multicolumn{2}{c|}{\smash{\raisebox{.5\normalbaselineskip}{$Z_4$}}}  \\  \hline 
\end{array},$$
where each $X,Y$ is a block of size $2\times 3$ and each $Z_r$ is a zero-sum block of size $r\times 4$.
Hence, we will construct $4(v+1)(t+1)$ blocks $X$, $2(t+1)$ blocks $Y$, $2v(w+1)(t+1)$  blocks $Z_4$ and $2(w+1)(t+1)$ blocks $Z_6$.

Take
$$\begin{array}{rcl}
\mathfrak{A} & = & \mathfrak{A}(8t+9, \; 8(3v+4)(t+1),\; 4(v+1)(t+1)),\\
\mathfrak{B} & = & \mathfrak{B}(4(6v+9)(t+1) ,\; t+1)   
\end{array}$$
as in Lemmas \ref{BlockA} and \ref{BlockB}, respectively.
Then, $\supp(\mathfrak{A})\dcup \supp(\mathfrak{B})=  [ 1, 16v(t+1)+24t+23 ]_2\dcup [ 16v(t+1)+24t+25, 32v(t+1)+48t+48 ]$, 
that is 
$$\supp(\mathfrak{A})\dcup \supp(\mathfrak{B}) = [1, mnc ] \setminus \left(\mathcal{M}_1\dcup\mathcal{M}_2  \dcup \mathcal{M}_{3} \right),$$
where
$$\begin{array}{rcl}
\mathcal{M}_1 & =& [ 2, 16v(t+1)+ 24t+24]_2, \\
\mathcal{M}_2 & =& [ 32v(t+1)+ 48t+49, (40v+16w)(t+1) + 48t+48  ],\\
\mathcal{M}_3 & =& [8(5v+2w)(t+1) + 48t+49 , 8(4vw+7v+6w) (t+1)+ 84t+84 ].
\end{array}$$
So, we replace the $4(v+1)(t+1)$ instances of $X$ with the arrays in  $\mathfrak{A}$ and the $2(t+1)$ instances of $Y$ with the arrays in $\mathfrak{B}$.
Hence, by construction, the first three columns have zero sum.

The set $\mathcal{M}_1\dcup \mathcal{M}_2$  can be written as a disjoint union 
 $F_1\dcup \ldots \dcup F_{\alpha}$, where $\alpha=(4v+4w+3)(t+1)$ and each $F_f$ is a $4$-set of type $2$.
The set $\mathcal{M}_3$ can be written as a disjoint union 
$G_1\dcup \ldots \dcup G_{\beta} \dcup H_1\dcup \ldots H_\gamma$, where $\beta=(8vw+4v+8w+7)(t+1)$,
$\gamma=t+1$, each $G_g$ is a $4$-set of type $1$ and each $H_h$ is a $8$-set of type $1$.

Construct $2(w+1)(t+1) $ blocks $Z_6$  as follows:
\begin{itemize}
\item[(1)] $(2w+1)(t+1)$ blocks $Q_1$ using the elements of two sets $F_f$ and four sets $G_g$;
\item[(2)] $t+1$ blocks $Q_2$ using the elements of one set $F_f$, three sets $G_g$ and one set $H_h$.
\end{itemize}
Construct $2v(w+1)(t+1) $ blocks $Z_4$  as follows:
\begin{itemize}
\item[(1)] $v(t+1)$ blocks $P_1$ using the elements of four sets $F_f$;
\item[(2)] $v(2w+1)(t+1)$ blocks $P_2$ using the elements of four sets $G_g$.
\end{itemize}
\end{proof}

\begin{ex}
We  aim to construct an $\IHS(10,7;2)$.
Following the proof of the previous proposition,
we first consider the two sets of $2\times 3$ zero-sum blocks $\mathfrak{A}(9, 56, 8)$ and $\mathfrak{B}(60,1)$. To construct the elements of these sets we use exactly the integers from $[ 1,39 ]_2 \dcup [41 , 80]$.
So, we still have to use the integers from
$\mathcal{M}_1=[2, 40]_2$, 
$\mathcal{M}_2=[81, 88]$ and
$\mathcal{M}_3=[89, 140]$.
We write $\mathcal{M}_1\dcup \mathcal{M}_2$ as
the disjoint union $F_1\dcup \ldots \dcup F_7$,
where $F_1=[2,8]_2$, $F_2=[10,16]_2$, 
$F_3=[18,24]_2$, $F_4=[26,32]_2$, $F_5=[34,40]_2$,
$F_6=[81,87]_2$ and $F_7=[82,88]_2$.
Finally, the set $\mathcal{M}_3$ can be written as
the disjoint union $G_1\dcup \ldots \dcup G_{11}\dcup H_1$, where
$G_1=[89,92]$, $G_2=[93,96]$,
$G_3=[97,100]$, $G_4=[101,104]$,
$G_5=[105,108]$, $G_6=[109,112]$,
$G_7=[113,116]$, $G_8=[125,128 ]$,
$G_9=[129,132 ]$, $G_{10}=[133,136]$,
$G_{11}=[137,140]$ and $H_1=[117,124] $.
We get 
 $$\begin{array}{|c|c|c|c|c|c|c|}\hline
1 & -61 & 60 & \sha 2 & \sha  -6 &  \sha -10 &  \sha 14 \\ \hline
-5 & 63 & -58 &  \sha -4 & \sha  8 & \sha  12 & \sha  -16 \\ \hline
-9 & 65 & -56 & \sha  -89 &  \sha 91 & \sha  93 &  \sha -95 \\ \hline
11 & -66 & 55 & \sha  90 & \sha  -92 & \sha  -94 & \sha  96 \\ \hline
-13 & 67 & -54 & \sha  -97 & \sha  99 & \sha  101 &  \sha -103 \\ \hline
15 & -68 & 53 &  \sha 98 &  \sha -100 & \sha  -102 &  \sha 104 \\ \hline
17 & -69 & 52 & 26 & -30 & -34 & 38 \\ \hline
-19 & 70 & -51 & -28 & 32 & 36 & -40 \\ \hline
-21 & 71 & -50 & -81 & 85 & 82 & -86 \\ \hline
23 & -72 & 49 & 83 & -87 & -84 & 88 \\ \hline
 \end{array},$$
  $$\begin{array}{|c|c|c|c|c|c|c|}\hline
3 & -62 & 59 &\sha  18 & \sha -20 & \sha -105 &\sha  107\\ \hline
-7 & 64 & -57 & \sha -22 & \sha 24 & \sha 106 &\sha  -108\\ \hline
-25 & 73 & -48 &\sha  -109 &\sha  110 &\sha  115 & \sha -116\\ \hline
27 & -74 & 47 &\sha  111 &\sha  -112 & \sha -113 & \sha 114\\ \hline
-29 & 75 & -46 &\sha\sha   124 &\sha  117 &\sha  -121 & \sha -120\\ \hline
31 & -76 & 45 &\sha  -122 & \sha -119 &\sha  118 &\sha  123\\ \hline
33 & -77 & 44 & 125 & -127 & -129 & 131\\ \hline
-35 & 78 & -43 & -126 & 128 & 130 & -132\\ \hline
-37 & 79 & -42 & -133 & 135 & 137 & -139\\ \hline
39 & -80 & 41 & 134 & -136 & -138 & 140\\ \hline
 \end{array}.$$
The two blocks $Z_6$ (highlighted in grey) are obtained by using 
$Q_1(F_1,F_2,G_1,G_2,G_3,G_4)$ and $Q_2(F_3,G_5, G_6 ,G_7,H_1)$;
the two blocks $Z_4$ are obtained by using 
 $P_1(F_4,F_5,F_6,F_7)$ and $P_2(G_8,G_9,G_{10},G_{11})$.
\end{ex}

\begin{prop}\label{n1}
Let $m,n,c$ be three positive integers such that $m\geqslant 6$ and $n\geqslant 9$.
Suppose that $m\equiv 2 \pmod 4$, $n\equiv 1\pmod 4$  and $c$ is even.
Then, there exists an $\IHS(m,n;c)$.
\end{prop}

\begin{proof}
Let $m=6+4v$, $n=9+4w$ and $c=2t+2$, where $v,w,t\geqslant 0$. 
Our $\IHS(m,n;c)$ will consist of $2t+2$ arrays of type
\begin{equation*}
\begin{array}{|c|ccc|cc|cc|cc|cccc}\hline
  Y      & \multicolumn{3}{c|}{} & \multicolumn{2}{c|}{} & \multicolumn{2}{c|}{} & \multicolumn{2}{c|}{} &\multicolumn{2}{c|}{}  \\ \cline{1-1}
  -X     & \multicolumn{3}{c|}{\smash{\raisebox{0\normalbaselineskip}{$W$}}}  & \multicolumn{2}{c|}{\smash{\raisebox{0\normalbaselineskip}{$Z_6$}}}  & \multicolumn{2}{c|}{\smash{\raisebox{0\normalbaselineskip}{$Z_6$}}} & \multicolumn{2}{c|}{\smash{\raisebox{0\normalbaselineskip}{$\cdots$}}} & \multicolumn{2}{c|}{\smash{\raisebox{0\normalbaselineskip}{$Z_6$}}}\\  \cline{1-1}
  -X      & \multicolumn{3}{c|}{}  & \multicolumn{2}{c|}{} & \multicolumn{2}{c|}{} & \multicolumn{2}{c|}{} &\multicolumn{2}{c|}{}   \\ \hline
   X     & \multicolumn{3}{c|}{}  &  \multicolumn{2}{c|}{} & \multicolumn{2}{c|}{} & \multicolumn{2}{c|}{} &\multicolumn{2}{c|}{}\\ \cline{1-1}
  -X     & \multicolumn{3}{c|}{\smash{\raisebox{.5\normalbaselineskip}{$(Z_6)^\T$}}}  &   \multicolumn{2}{c|}{\smash{\raisebox{.5\normalbaselineskip}{$Z_4$}}}   &   \multicolumn{2}{c|}{\smash{\raisebox{.5\normalbaselineskip}{$Z_4$}}}   &   \multicolumn{2}{c|}{\smash{\raisebox{.5\normalbaselineskip}{$\cdots$}}}   &   \multicolumn{2}{c|}{\smash{\raisebox{.5\normalbaselineskip}{$Z_4$}}}  \\  \hline 
  \vdots    & \multicolumn{3}{c|}{\vdots}   &     \multicolumn{2}{c|}{\vdots} & \multicolumn{2}{c|}{\vdots} & \multicolumn{2}{c|}{\ddots} &\multicolumn{2}{c|}{\vdots}  \\[1mm] \hline
   X   & \multicolumn{3}{c|}{}    &  \multicolumn{2}{c|}{} & \multicolumn{2}{c|}{} & \multicolumn{2}{c|}{} &\multicolumn{2}{c|}{}  \\ \cline{1-1}
  -X    &\multicolumn{3}{c|}{\smash{\raisebox{.5\normalbaselineskip}{$(Z_6)^\T$}}} &   \multicolumn{2}{c|}{\smash{\raisebox{.5\normalbaselineskip}{$Z_4$}}}   &   \multicolumn{2}{c|}{\smash{\raisebox{.5\normalbaselineskip}{$Z_4$}}}   &   \multicolumn{2}{c|}{\smash{\raisebox{.5\normalbaselineskip}{$\cdots$}}}   &   \multicolumn{2}{c|}{\smash{\raisebox{.5\normalbaselineskip}{$Z_4$}}}  \\  \hline 
\end{array},
\end{equation*}
where each $X,Y$ is a block of size $2\times 3$, each $Z_r$ is a zero-sum block of size $r\times 4$ and each $W$ is zero-sum block of size $6\times 6$.
Hence, we will construct $4(v+1)(t+1)$ blocks $X$, $2(t+1)$ blocks $Y$, 
$2vw(t+1)$  blocks $Z_4$,  $2(v+w)(t+1) $ blocks $Z_6$ and $2(t+1)$ blocks $W$.

Take
$$\begin{array}{rcl}
\mathfrak{A} & = & \mathfrak{A}(8t+9, \; 8(3v+4)(t+1),\; 4(v+1)(t+1)),\\
\mathfrak{B} & = & \mathfrak{B}(4(6v+9)(t+1) ,\; t+1)   
\end{array}$$
as in Lemmas \ref{BlockA} and \ref{BlockB}, respectively.
Then, 
$$\supp(\mathfrak{A})\dcup \supp(\mathfrak{B})= [1, mnc ]
\setminus \left(\mathcal{M}_1\dcup\mathcal{M}_1 \dcup\mathcal{M}_3  \right),$$
where
$$\begin{array}{rcl}
\mathcal{M}_1 & =& [ 2, 16v(t+1)+ 24t+24]_2, \\
\mathcal{M}_2 & =& [ 32v(t+1)+ 48t+49, 32v(t+1)+72t+72 ],\\
\mathcal{M}_3 & =& [32v(t+1)+72t+73,  8(4vw+9v+6w)(t+1) + 108t+108 ].\\
\end{array}$$
So, we replace the $4(v+1)(t+1)$ instances of $X$ with the arrays in  $\mathfrak{A}$ and the $2(t+1)$ instances of $Y$ with the arrays in $\mathfrak{B}$.
By construction, the first three columns have zero sum.

The set $\mathcal{M}_1\dcup \mathcal{M}_2$  can be written as a disjoint union 
$F_1\dcup \ldots \dcup F_{\alpha}\dcup J_1\dcup \ldots \dcup J_{\delta}$,
where $\alpha=(2v+3)(t+1)$, 
$\delta=3(t+1)$, each $F_f$ is a $4$-set of type $2$ and each $J_j$ is a $8$-set of type $2$.
The set $\mathcal{M}_3$ can be written as a disjoint union 
$G_1\dcup \ldots \dcup G_\beta \dcup H_1\dcup \ldots H_\gamma$, 
where $\beta=(8vw+10v+12w+3  )(t+1)$, 
$\gamma=3(t+1)$, each $G_g$ is a $4$-set of type $1$ and each $H_h$ is a $8$-set of type $1$.

Construct $2(t+1)$ blocks $W$  as follows:
\begin{itemize}
\item[(1)] $t+1$ blocks $R_1$ using the elements of three sets $F_f$ and three sets $J_j$;
\item[(2)] $t+1$ blocks $R_2$ using the elements of three sets $G_g$ and  three sets $H_h$.
\end{itemize}
Construct $2(v+w)(t+1)$ blocks $Z_6$ as follows:
\begin{itemize}
\item[(1)] $v(t+1)$ blocks $Q_1$ using the elements of two sets $F_f$ and four sets $G_g$;
\item[(2)] $(v+2w)(t+1)$ blocks $Q_3$ using the elements of six sets $G_g$.
\end{itemize}
Finally, construct $2vw(t+1)$ blocks $Z_4$  taking blocks $P_2$ and using, for each of them, the elements of four sets~$G_g$.
\end{proof}

\begin{proof}[Proof of Theorem 
\ref{mainB}]
$(1)$ If $m$ and $n$ are two even integers, then $mn\equiv 0 \pmod 4$. By \cite[Lemma 1.3]{ABD}  there exists a shiftable 
$\H(m,n)$, say $A$.
The arrays $A,A\pm mn, \ldots, A\pm (c-1)mn$ 
constitute  an $\IHS(m,n;c)$.

\noindent $(2)$ The existence of an $\IHS(m,n;c)$ when $m\equiv 0 \pmod 4$ and $n $ is an odd integer was proved in \cite{MP3}, except for $n=5$.
The existence of an $\IHS(m,5;c)$ with $m\equiv 0 \pmod 4$ follows from Lemma~\ref{45even} and Lemma~\ref{45odd}.

\noindent $(3)$ Note that the condition $mnc\equiv 0,3\pmod 4$ implies that $c$ is an even integer, as $m\equiv 2 \pmod 4$ and $n$ is odd.
Hence, this case follows from Propositions \ref{n3} and \ref{n1}.

\noindent $(4)$ This case follows from \cite[Theorem 1.8]{PT}.
\end{proof}

\begin{proof}[Proof of Theorem 
\ref{mainC}]
By Theorem \ref{diago} and
Theorem \ref{mainA}, there exists an integer $\H(m,n;s,k)$ whenever 
$\gcd(s,k)>1$, proving case $(1)$.
So, without loss of generality, we can assume that $k$ is an odd integer such that $\gcd(s,k)=1$.
From $ms=nk$ we obtain
$m=ck$ and $n=cs$ for some $c\geqslant 1$.
It is clear that, to construct an 
integer $\H(ck,cs;s,k)$, it suffices to use the blocks of an 
$\IHS(k,s;c)$.
Hence, case $(2)$ follows from item $(2)$ of Theorem \ref{mainB}, while case $(3)$ follows from  items $(3)$ and $(4)$ of the same theorem.
\end{proof}

\section*{Acknowledgments}

Both authors are partially supported by INdAM-GNSAGA.


\begin{thebibliography}{20}

\bibitem{A} D.S. Archdeacon,
Heffter arrays and biembedding graphs on surfaces,
\textit{Electron. J. Combin.} \textbf{22} (2015), \#P1.74.

\bibitem{ABD} D.S. Archdeacon, T. Boothby, J.H. Dinitz,
Tight Heffter arrays exist for all possible values,
\textit{J. Combin. Des.} \textbf{25} (2017), 5--35.

\bibitem{ADDY} D.S. Archdeacon, J.H. Dinitz, D.M. Donovan, E.S. Yaz{\i}c{\i},
Square integer Heffter arrays with empty cells,
\textit{Des. Codes Cryptogr.} \textbf{77} (2015), 409--426.

\bibitem{C} S. Cichacz,     
A $\Gamma$-magic rectangle set and group distance magic labeling,
in: \textit{Combinatorial algorithms, IWOCA 2014} (Eds. J. Kratochvíl, M. Miller,  D. Froncek),
Lecture Notes in Comput. Sci., 8986 Springer, Cham, 2015.

\bibitem{DW} J.H. Dinitz, I.M. Wanless,
The existence of square integer Heffter arrays,
\textit{Ars Math. Contemp.} \textbf{13} (2017), 81--93. 

\bibitem{F13}  D. Froncek, 
Handicap distance antimagic graphs and incomplete tournaments,
\textit{AKCE Int. J. Graphs Comb.} \textbf{10} (2013), 119--127.

\bibitem{MP} F. Morini, M.A. Pellegrini,
On the existence of integer relative Heffter arrays,
\textit{Discrete Math.} \textbf{343} (2020), \#112088.

\bibitem{MP3} F. Morini, M.A. Pellegrini,
Rectangular Heffter arrays: a reduction theorem,
\textit{Discrete Math.} \textbf{345} (2022), \#113073.

\bibitem{MP4} F. Morini, M.A. Pellegrini,
Magic partially filled arrays on abelian groups,
\textit{J. Combin. Des.} \textbf{31} (2023), 347--367.

\bibitem{PT} M.A. Pellegrini, T. Traetta,
Towards a solution of Archdeacon's conjecture on integer Heffter arrays, 
\textit{J. Combin. Des.} \textbf{33} (2025), 310--323. 


\end{thebibliography}
\end{document}